\documentclass[11pt]{article}
\usepackage{amstext,amssymb,amsmath,amsbsy}

\usepackage{hyperref}
\usepackage{amscd}
\usepackage{amsfonts}
\usepackage{indentfirst}
\usepackage{verbatim}
\usepackage{amsmath}
\usepackage{amsthm}
\usepackage{enumerate}
\usepackage{graphicx}
\usepackage{subfig}
\usepackage{color}
\usepackage[OT1]{fontenc}
\usepackage[latin1]{inputenc}
\usepackage[english]{babel}
\usepackage{amssymb}
\newtheorem{theorem}{Theorem}
\newtheorem{lemma}{Lemma}

\newtheorem{proposition}{Proposition}

\newtheorem{remark}{Remark}

\newtheorem{definition}{Definition}
\numberwithin{equation}{section}
\numberwithin{theorem}{section}
\numberwithin{lemma}{section}
\numberwithin{notation}{section}
\numberwithin{proposition}{section}
\numberwithin{corollary}{section}
\numberwithin{corollary}{section}
\numberwithin{example}{section}
\numberwithin{definition}{section}
\numberwithin{remark}{section}

\newcommand{\proofend}{\hfill $\Box$ }

\newcommand{\sle}{\lesssim}

\newcommand{\cV}{ \hat   V}

\newcommand{\supp}{\operatorname{supp}}

\newcommand{\dive}{\operatorname{div}}

\newcommand{\eps}{\varepsilon}

\newcommand{\epss}{s}

\newcommand{\loc}{_{loc}}

\newcommand{\mC}{\mathbb{C}}

\newcommand{\mR}{\mathbb{R}}

\newcommand{\mc}{\mathrm{c}}

\newcommand{\hu}{\hat u}

\newcommand{\cA}{{\hat A}}
\newcommand{\cS}{\hat \Sigma}
\newcommand{\cU}{\hat u}

\title{Cloaking via anomalous localized resonance for doubly complementary media in the finite frequency regime}

\author{Hoai-Minh Nguyen \footnote{EPFL SB MATHAA CAMA, Station 8,  CH-1015 Lausanne, hoai-minh.nguyen@epfl.ch}}

\begin{document}

\maketitle

\begin{abstract} Cloaking  a source via anomalous localized resonance (ALR)  was discovered by Milton and Nicorovici in \cite{MiltonNicorovici}. 
A general setting in which cloaking a source via ALR takes place
is the settting of doubly complementary media. This was  introduced and studied in \cite{Ng-CALR}  for the quasistatic regime. 
In this paper, we study cloaking a source via ALR for doubly complementary media in  the finite frequency regime as a natural continuation of \cite{Ng-CALR}. We establish the following results: 
1) Cloaking a source via ALR appears if and only if the power blows up;  2) The power blows up if the source is ``placed" near the plasmonic structure;  3) The power remains bounded if the source is far away from the plasmonic structure. 
Concerning the analysis,  we extend  ideas from \cite{Ng-CALR} and  add new insights on the problem which allows us to overcome difficulties related to the finite frequency regime and  to
obtain new information on the problem. In particular, we are able to  characterize the behaviour of the fields far enough from the plasmonic shell as the loss goes to 0 for an {\bf arbitrary source} outside the core-shell structure in the doubly complementary media setting. 
\end{abstract}


\section{Introduction}
Negative index materials (NIMs) were first investigated theoretically by Veselago in \cite{Veselago}. The  existence of such materials was confirmed by Shelby, Smith, and Schultz in \cite{ShelbySmithSchultz}. The study of NIMs has attracted a lot attention in the scientific community thanks to their many applications. One of the appealing ones is cloaking. There are at least three ways to do cloaking using NIMs. The first one is based on plasmonic structures introduced by Alu and Engheta in \cite{AluEngheta}. The second one uses the concept of complementary media. This was suggested by Lai et al. in \cite{LaiChenZhangChanComplementary} and confirmed theoretically  in  \cite{Ng-Negative-Cloaking} for inspired schemes. 
The last one  is based on the concept of ALR discovered by Milton and Nicorovici in \cite{MiltonNicorovici}. In this paper, we concentrate on the last method.  

\medskip
Cloaking a source  via ALR  was discovered by Milton and 
 Nicorovici in \cite{MiltonNicorovici}. Their work has root from \cite{NicoroviciMcPhedranMilton94} (see also \cite{NicoroviciMcPhedranMiltonPodolskiy1}) where the localized resonance was observed and established for constant symmetric plasmonic structures  in the two dimensional quasistatic regime. 
More precisely, in \cite{MiltonNicorovici}, the authors studied  core-shell plasmonic structures  in which  a circular shell has permittivity $-1 - i \delta$ while its complement has permittivity 1  where  $\delta$ denotes the loss of the material in the shell \footnote{In fact, in \cite{MiltonNicorovici} and in other works, the authors consider the permittivity $-1 + i \delta$ instead of $-1 - i \delta$; but this point is not essential.}.  Let $r_1$ and $r_2$ be the inner and the outer radius of the shell. 
They showed that there is a critical radius  $r_* := (r_2^3 r_1^{-1})^{1/2}$ such that  
a dipole is not seen by an observer away from the core-shell structure, hence it is cloaked,  if and only if the dipole is within distance $r_*$ of the shell. Moreover, 
the power $E_\delta(u_\delta)$  of the field $u_\delta$, which is defined in  \eqref{def-power},  blows up as the loss $\delta$ goes to 0.  Two key features of this phenomenon are: 1)  the localized resonance, i.e., the fields blow up in some regions and remain bounded in some others as the loss goes to 0; 2) the connection between the localized resonance and the blow up of the power as the loss goes to 0.

\medskip

Cloaking a source via ALR has been mainly studied  in the quasistatic regime.  In  \cite{BouchitteSchweizer10}, Bouchitte and Schweizer  proved that a small circular inclusion of radius $\gamma(\delta)$ (with $\gamma(\delta) \to 0$ fast enough) is cloaked by the core-shell plasmonic structure  mentioned above in the two dimensional quasistatic regime  if the inclusion is located within distance $r_*$ of the shell. Otherwise it is visible. 
Concerning the second feature of cloaking a source via ALR, the blow up of the power  was studied for a more general setting in two dimensional quasistatic regime by Ammari et al. in \cite{AmmariCiraoloKangLeeMilton} and  Kohn et al. in \cite{KohnLu}. More precisely, they considered  non-radial core-shell structures in which the shell has permittivity $-1 - i \delta$ and its complement has permitivity 1. In \cite{AmmariCiraoloKangLeeMilton},  Ammari et al. dealt with   arbitrary shells  and provided a characterization  of sources for which the power blows up via the information of the spectral decomposition of the Neumann-Poincar\'e type operator.  In \cite{KohnLu},  Kohn et al. considered core-shell structures  in which the outer boundary of the shell is round  but the inner is not and  established the blow up of the power for some class of sources using a variational approach.  A connection between the blow up of the power and the localized resonance subtly depends  on the geometry and property of plasmonic structures. It was showed in \cite{MinhLoc1} that such a connection does not hold in general. Cloaking a source via ALR in some special three dimensional  geometry was studied in \cite{AmmariCiraoloKangLeeMilton2}. Motivated by the concept of reflecting complementary media suggested and studied in \cite{Ng-Complementary} and results mentioned above, in  \cite{Ng-CALR} we studied cloaking a source via ALR for a general core shell structure of  doubly complementary media property (see Definition~\ref{def-DCM}) in the quasistatic regime. \footnote{Roughly speaking, the plasmonic shell is not only complementary with a part of the complement of the core shell but also complement to a part of the core.}  More precisely,  we  established the following three properties for doubly complementary media: 
\begin{itemize}
\item[P1)] Cloaking a source via ALR appears if and only if the power blows up. 
\item[P2)] The power blows up if the source is located  ``near" the shell. 
\item[P3)] The power remains bounded if the source is far away from the shell. 
\end{itemize}
Using these results, we extended various results mentioned previously. Moreover, 
we were able to obtain  schemes  to cloak  an arbitrary source concentrating on an arbitrary smooth bounded manifold of codimension 1 placed in an   arbitrary medium via ALR; the cloak is independent of the source. The analysis in \cite{Ng-CALR} is on one hand based on the reflecting techniques initiated in \cite{Ng-Complementary}, the removing localized technique introduced in \cite{Ng-Superlensing, Ng-Negative-Cloaking} to deal with the localized resonance. On the other hand, the analysis in \cite{Ng-CALR} is based on  new observations on the Cauchy problems and the separation of variables technique for a general shell introduced there. The implement of this technique is an ad-hoc part of \cite{Ng-CALR}. 

\medskip 
In this paper, we study cloaking a source via ALR for the finite frequency regime as a natural continuation of \cite{Ng-CALR}. More precisely, we establish Properties P1), P2) and P3)  for doubly complementary media in the finite frequency regime. As a consequence, we are also able, as in the quasistatic regime,  to obtain  schemes to cloak  {\bf a non generic  arbitrary source} concentrating on an {\bf arbitrary smooth bounded manifold of codimension 1} placed in an {\bf  arbitrary medium} via ALR; the cloak is independent of the source (see Section~\ref{sect-cloaking}).
Concerning the analysis,  we extend ideas from \cite{Ng-CALR} and add various new insights for the problem which allows us to overcome difficulties related to the finite frequency regime such as the use of the maximum priniciple, to shorten the approach in \cite{Ng-CALR}, and to
obtain new information on the cloaking a source via ALR. In particular, we can characterize the behaviour of the fields far enough from the plasmonic shell as the loss goes to 0 for arbitrary sources in the doubly complementary media setting (Theorem~\ref{thm-main}). This fact  is interesting in itself and new to our knowledge.

\medskip 
 Let $k > 0$ and let $A$ be a (real) uniformly elliptic symmetric matrix defined in $\mR^d$ ($d \ge 2$), and $\Sigma$ be a  bounded real function defined in $\mR^d$.  Assume that  
\begin{equation}\label{cond-I}
A(x) = I  \quad \Sigma(x) = 1  \mbox{ for $|x|$ large enough} 
\end{equation}
and 
\begin{equation}\label{cond-C1}
A \mbox{ is piecewise } C^1.   
\end{equation}
Let $\Omega_1 \subset \subset \Omega_2 \subset \subset \mR^d$ be smooth connected open subsets of $\mR^d$,  and set, for $\delta \ge 0$,  
\begin{equation}\label{def-sd}
s_\delta (x) = \left\{\begin{array}{cl} - 1 - i \delta & \mbox{ in } \Omega_2 \setminus \Omega_1, \\[6pt]
1 & \mbox{ in } \mR^d \setminus (\Omega_2 \setminus \Omega_1). 
\end{array}\right.
\end{equation}
For $f \in L^2_{\mc}(\mR^d)$ with   $\supp f \cap B_{r_2} = \O$ and $\delta > 0$,  let  
$u_\delta \in H^1_{\loc}(\Omega)$ be the unique outgoing  solution to 
\begin{equation}\label{def-ud}
\dive (s_\delta A \nabla u_\delta) + k^2 s_0 \Sigma u_\delta = f \mbox{ in } \mR^d. 
\end{equation}
Here and in what follows 
\begin{equation*}
L^2_{\mc}(\mR^d) : = \big\{ f \in L^2(\mR^d) \mbox{ with compact support} \big\}. 
\end{equation*}
For $R>0$ and $x \in \mR^d$, we also denote $B(x, R)$  the ball in $\mR^d$ centered at $x$ and  of radius $R$; when $x = 0$, we simply denote $B(x, R) $ by $B_R$. Recall that a function $u \in H^1_{\loc}(\mR^d \setminus B_R)$ for some $R>0$ which is a solution to the equation $\Delta u + k^2 u = 0 $ in $\mR^d \setminus B_R$ is said to satisfy the outgoing condition if 
\begin{equation*}
\partial_r u - i k u = o( r^{\frac{1-d}{2}}) \mbox{ as } r = |x| \to + \infty. 
\end{equation*}
The power $E_\delta(u_\delta)$ is defined by (see, e.g., \cite{MiltonNicorovici})
\begin{equation}\label{def-power}
E_\delta (u_\delta) = \delta  \int_{\Omega_2 \setminus \Omega_1} |\nabla u_\delta|^2. 
\end{equation}
The normalization of  $u_\delta$ is  $v_\delta \in H^1_{\loc}(\mR^d)$ which is  the unique outgoing solution of 
\begin{equation}\label{def-vd}
\dive (s_\delta A \nabla v_\delta) + k^2 s_0 \Sigma v_\delta = f_\delta \mbox{ in } \mR^d. 
\end{equation}
Here 
\begin{equation*}
f_\delta = c_\delta f, 
\end{equation*}
 and  $c_\delta$ is the normalization  constant such that 
\begin{equation}\label{Power-vn}
\delta^{1/2} \int_{B_{r_2} \setminus B_{r_1}} |\nabla v_\delta|^2 = 1. 
\end{equation}

In this paper, we establish properties P1), P2), and P3) for $(A, \Sigma)$ of doubly complementary property.  Before giving the definition of doubly complementary media for a general core-shell structure, let us recall the definition of reflecting complementary media introduced  in \cite[Definition 1]{Ng-Complementary}. 

\begin{definition}[Reflecting complementary media]   \fontfamily{m} \selectfont
 \label{def-Geo} Let $\Omega_1 \subset \subset \Omega_2 \subset \subset \Omega_3 \subset \subset \mR^d$ be smooth open bounded subsets of $\mR^d$. The media $(A, \Sigma)$ in $\Omega_3 \setminus \Omega_2$ and $(-A, -\Sigma)$ in $\Omega_2 \setminus \Omega_1$  are said to be  {\it reflecting complementary} if 
there exists a diffeomorphism $F: \Omega_2 \setminus \bar \Omega_1 \to \Omega_3 \setminus \bar \Omega_2$ such that 
\begin{equation}\label{cond-ASigma}
(F_*A, F_*\Sigma)   = (A, \Sigma )   \mbox{ for  } x \in  \Omega_3 \setminus \Omega_2, 
\end{equation}
\begin{equation}\label{cond-F-boundary}
F(x) = x \mbox{ on } \partial \Omega_2, 
\end{equation}
and the following two conditions hold: 1) There exists an diffeomorphism extension of $F$, which is still denoted by  $F$, from $\Omega_2 \setminus \{x_1\} \to \mR^d \setminus \bar \Omega_2$ for some $x_1 \in \Omega_1$; 2) There exists a diffeomorphism $G: \mR^d \setminus \bar \Omega_3 \to \Omega_3 \setminus \{x_1\}$ such that 
\begin{equation}\label{cond-G-boundary}
\quad G(x) = x \mbox{ on } \partial \Omega_3,
\end{equation}
and
\begin{equation}\label{extension}
G \circ F : \Omega_1  \to \Omega_3 \mbox{ is a diffeomorphism if one sets } G\circ F(x_1) = x_1.
\end{equation}
\end{definition}

Here and in what follows, if $T$ is a diffeomorphism and  $a$ and $\sigma$ are a matrix-valued function and a complex function, we use the following standard notations  
\begin{equation}\label{def-F*}
T_*a(y) = \frac{DT(x) a(x) DT(x)^T}{|\det DT(x)|} \quad \mbox{ and } \quad T_*\sigma(y) = \frac{ \sigma(x) }{|\det DT(x)|} \quad  \mbox{ where } x = T^{-1}(y). 
\end{equation}

In \eqref{cond-F-boundary} and \eqref{cond-G-boundary}, $F$ and $G$ denote some  diffeomorphism extensions of $F$ and $G$ in a neighborhood of $\partial \Omega_2$ and of $\partial \Omega_3$. Conditions  \eqref{cond-ASigma} and \eqref{cond-F-boundary} are the main assumptions in Definition~\ref{def-Geo}. The key point behind this requirement is roughly speaking the following property:  if $u_0 \in H^1(\Omega_3 \setminus \Omega_1)$ is a solution of $\dive(s_0 A \nabla u_0) + k^2 s_0 \Sigma u_0= 0$ in $\Omega_3 \setminus \Omega_2$ and if $u_1$ defined in $\Omega_3 \setminus \Omega_2$ by  $u_1 = u_0 \circ F^{-1}$, then $\dive(A \nabla u_1) + k^2 \Sigma u_1 = 0$ in $\Omega_3 \setminus \Omega_2$, $u_1 - u_0 = A \nabla (u_1 - u_0) \cdot \nu  = 0 $ on $\partial \Omega_2$ (see~\cite[Lemma 2]{Ng-Complementary}). Hence $u_1 = u$ in $\Omega_3 \setminus \Omega_2$ by the unique continuation principle.   Conditions 1) and 2)  are  mild assumptions. Introducing $G$ makes the analysis more accessible, see \cite{Ng-Complementary, Ng-Superlensing, Ng-Negative-Cloaking, MinhLoc2} and the analysis presented in this paper. 

\begin{remark} \fontfamily{m} \selectfont The class of reflecting complementary media has played an important role in the other applications of NIMs such as cloaking and superlensing using complementary see \cite{Ng-Superlensing, Ng-Negative-Cloaking, MinhLoc2}. 
\end{remark}

\begin{remark} \fontfamily{m} \selectfont  Let $d=2$,  $A = I$, $0 < r_1 < r_2 < + \infty$ and set $r_3 = r_2^2/r_1$. Letting $F$ be the Kelvin transform with respect to $\partial B_{r_2}$, i.e., $F(x) = r_2^2 x/ |x|^2$ and $\Omega_i = B_{r_i}$,  one can  verify that in the quasistatic regime  the core-shell structures considered by Milton and Nicorovici in \cite{MiltonNicorovici} and by Kohn et al. in \cite{KohnLu}  have the reflecting complementary property.  

\end{remark}

We are ready to introduce the concept of  doubly complementary media for the finite frequency regime. 

\begin{figure}[h!]
\begin{center}
\includegraphics[width=6cm]{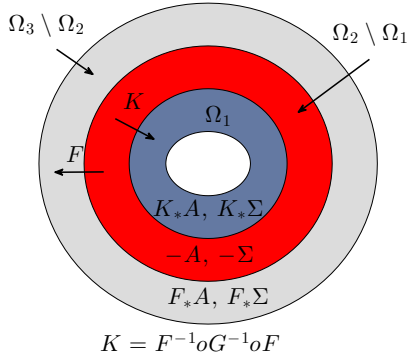} 
\caption{$(s_0A, s_0\Sigma)$ is doubly complementary: $(-A, -\Sigma)$ in $\Omega_2 \setminus \Omega_1$ (the red region) is complementary to $(F_*A, F_*\Sigma)$ in $\Omega_3 \setminus \Omega_2$ (the grey region) and $(K_*A, K_*\Sigma)$ with $K = F^{-1} \circ G^{-1} \circ F$ in $K(B_{r_2} \setminus B_{r_1})$ (the blue grey region).} \label{fig1}
\end{center}
\end{figure}

\begin{definition} \label{def-DCM} \fontfamily{m} \selectfont The medium $(s_0 A, s_0\Sigma)$ is said to be {\it doubly complementary} if  for some $\Omega_2 \subset \subset \Omega_3$, $(A, \Sigma)$ in $\Omega_3 \setminus \Omega_2$ and $(-A, - \Sigma)$ in $\Omega_2 \setminus \Omega_1 $ are reflecting complementary, and 
\begin{equation}\label{def-DC}
F_*A = G_* F_*A = A  \quad \mbox{ and } \quad F_*\Sigma = G_* F_*\Sigma = \Sigma \mbox{ in } B_{r_3} \setminus B_{r_2}, 
\end{equation}
for some $F$ and $G$ coming from Definition~\ref{def-Geo} (see Figure \ref{fig1}). 
\end{definition}

The reason for which media satisfying \eqref{def-DC} are called doubly complementary media is that $(-A, -\Sigma)$ in $B_{r_2} \setminus B_{r_1}$ is not only complementary to $(A, \Sigma)$ in $\Omega_3 \setminus \Omega_2$ but also to $(A, \Sigma)$ in $(G \circ F)^{-1}(\Omega_3 \setminus \overline{\Omega_2})$ (a subset of $\Omega_1$) (see \cite{Ng-CALR-CRAS}). The key property behind Definition~\ref{def-DCM} is as follows. Assume $u_0 \in H^1_{\loc}(\mR^d)$ is a solution of \eqref{def-ud} with $\delta = 0$ and  set $u_1 = u \circ F^{-1}$ and $u_{2} = u_1 \circ G^{-1}$. Then $u_1, u_2$ satisfy the  equation 
$\dive(A \nabla \cdot) + k^2 \Sigma \cdot = 0$  in $\Omega_3 \setminus \Omega_2$ (see \cite[Lemma 2]{Ng-Complementary}). 

\begin{remark}\fontfamily{m} \selectfont
Taking $d=2$, $A = I$ and $r_3 = r_2^2/r_1$, and letting $F$ and $G$ be the Kelvin transform with respect to $\partial B_{r_2}$ and $\partial B_{r_3}$,  one can  verify that the core-shell structures considered by Milton and Nicorovici in \cite{MiltonNicorovici} is of doubly complementary property.  It is worthy to note that one requires no information of $A$ outside $B_{r_3}$ and inside $B_{r_1^2/ r_2}$ in the definition of doubly complementary media. 
\end{remark}

\begin{remark} \label{rem-DCM} \fontfamily{m} \selectfont Given $(A, \Sigma)$ in $\mR^d$ it is not easy in general to verify whether or not $(s_0A, s_0 \Sigma)$ is doubly complementary. Nevertheless, given $\Omega_1 \subset \Omega_2 \subset \subset \Omega$ and  $(A, \Sigma)$ in $\Omega_3 \setminus \Omega_2$, it is quite easy to choose $(A, \Sigma)$ in $\Omega_2$ such that $(s_0A, s_0\Sigma)$ is doubly complementary. One just needs to  choose diffeomorphisms $F$ and $G$ as in Definition~\ref{def-Geo} and define $(A, \Sigma) = (F^{-1}_*A, F^{-1}_*\Sigma)$ in $\Omega_2 \setminus \Omega_1$ and $(A, \Sigma) = (G^{-1}_*F^{-1}_*A, G^{-1}*F^{-1}_*\Sigma)$ in $G^{-1}\circ F^{-1} (\Omega_3 \setminus \Omega_2)$. This idea is used  in Section~\ref{sect-cloaking} when we discuss cloaking sources in an arbitrary medium. 

\end{remark}

The main result of this paper is the following theorem which reveals the behavior of $u_\delta$ for a general source $f$.

\begin{theorem}\label{thm-main}
Let $d \ge 2$, $k>0$, $f \in L^2_{\mc}(\mR^d)$ with $\supp f \cap \Omega_2 = \O$,   and let $u_\delta \in H^1_{\loc}(\mR^d)$ be the unique outgoing solution of \eqref{def-ud}.  Assume that $(s_0A, s_0\Sigma)$ is doubly complementary. 
Then 
\begin{equation}\label{part2}
u_\delta \to \cU\mbox{ weakly in } H^1_{\loc}(\mR^d \setminus B_{r_3}), 
\end{equation}
where $\cU \in H^1_{\loc}(\mR^d)$ is the unique outgoing solution of 
\begin{equation}\label{UU}
\dive (\cA \nabla {\cU}) + k^2 \cS \cU= f \mbox{ in } \mR^d. 
\end{equation}
Here 
\begin{equation}\label{def-cAS}
(\cA, \cS) : = \left\{\begin{array}{cl} (A, \Sigma) & \mbox{ in } \mR^d \setminus \Omega_3, \\[6pt]
(G_*F_*A, G_*F_*\Sigma) & \mbox{ in } \Omega_3. 
\end{array}\right.
\end{equation}
\end{theorem}

Using Theorem~\ref{thm-main}, one can establish the equivalence between the blow up of the power and the cloaking via ALR as follows. Suppose that the power blows up,  i.e.,  
\begin{equation*}
\lim_{n \to \infty} \delta_n \| \nabla u_{\delta_n}\|_{L^2(B_{r_2} \setminus B_{r_1})}^2 = + \infty. 
\end{equation*}
Then, by Theorem~\ref{thm-main},  $v_{\delta_n} \to 0$ in $\mR^d \setminus B_{r_3}$; the localized resonance takes place. The source $\alpha_{\delta_n} f$ is not seen by observers far away from the shell: the source is cloaked. If the power $E_{\delta_n} (u_{\delta_n})$ remains bounded,  then $u_{\delta_n} \to \cU$ weakly in $H^1_{\loc}(\mR^d \setminus B_{r_3})$. Since $\cU \in H^1_{\loc}(\mR^d)$ is the unique outgoing solution to \eqref{UU}, the source is not cloaked.  

\begin{remark} \fontfamily{m} \selectfont
It follows from \eqref{def-F*} that if  $(s_0A, s_0\Sigma)$ is doubly complementary media then $(A, \Sigma)$ is not piecewise constant. This is the reason for which there are very few results on cloaking via ALR in the finite frequency regime  in both two and three dimensions. 
\end{remark}

In comparison with results in \cite{Ng-CALR}, Theorem~\ref{thm-main} is stronger: no conditions on the blow up rate of the power is required. The proof of Theorem~\ref{thm-main} is in the spirit of \cite{Ng-CALR}. Nevertheless, we incorporate two important ingredients. The first one is on the blow up rate of the power of $u_\delta$ in \eqref{est5-TM} which is derived in this paper but was assumed previously. 
The second one is on the removing localized singularity technique.  Here, we are able to construct in a simple and  robust way the singular part which is necessary to removed. This  helps us to avoid the ad-hoc separation of variables for a general shell developed and implemented in \cite{Ng-CALR}. 
The construction of the removing term comes from a remark of Etienne Sandier. The author would like to thank him for it.  To our knowledge, Theorem~\ref{thm-main} is new and is the first result providing the connection between the blow up of the power and the invisibility of a source in the finite frequency regime.

\medskip

Concerning the blow up of the power, we can prove the following result which holds for a large class of media in which the reflecting complementary property holds only locally. 

\begin{proposition}\label{pro1} Assume that there exists  a diffeomorphism $F: \Omega_2 \setminus \Omega_1 \to \Omega_3 \setminus \Omega_2$ for some $\Omega_2 \subset \subset \Omega_3 \subset \mR^d$ such that $F(x) = x$ on $\partial \Omega_2$, 
\begin{equation*}
(A, \Sigma) = (F_*A, F_*\Sigma) \mbox{ in } D \mbox{ where } D : = B(x_0, R_0) \cap (\Omega_3 \setminus \Omega_2). 
\end{equation*}
for some $x_0 \in \partial \Omega_2$ and $R_0> 0$. Let $f \in L^2_{\mc}(\mR^d)$ and assume that $A$ is Lipschitz in $\overline D$. There exists $0 < r_0 < R_0$, independent of $f$,  such that if  there is no solution $v \in H^1(D_1)$ where $D_1 := D \cap B(x_0, r_0)$ to the Cauchy problem: 
\begin{equation*}
\dive(A \nabla v) + k^2 \Sigma v = f \mbox{ in } D_1  \quad \mbox{ and } \quad v = A \nabla v \cdot \nu = 0 \mbox{ on } \partial D_1\setminus  \partial B(x_0, r_0),  
\end{equation*}
then   
\begin{equation*}
\limsup_{\delta \to 0} \delta \int_{\Omega_2 \setminus \Omega_1} |\nabla u_\delta|^2 = + \infty, 
\end{equation*}
where $u_\delta \in H^1_{\loc}(\mR^d)$ is the unique outgoing solution of \eqref{def-ud}. 
\end{proposition}

The proof of Proposition~\ref{pro1} essentially uses the ideas in the one of  \cite[Lemma 10]{Ng-WP} which has root from \cite{Ng-CALR}. It is clear that Proposition~\ref{pro1} implies Property P2). Note that $(s_0A, s_0\Sigma)$ is not required to be doubly complementary in Proposition~\ref{pro1}. 

\medskip 
Concerning the boundedness of the power, we have the following result, which implies Property P3). 

\begin{proposition} \label{pro2} Let $d=2, 3$, $0< \delta < 1$,  $f \in L^2_{\mc}(\mR^d)$, and   let $u_\delta \in H^1(\mR^d)$ be the unique solution \eqref{def-ud}. 
Assume that $(s_0A, s_0\Sigma)$ is doubly complementary media. We have, if $\supp f \cap \Omega_3 = \O$ then 
\begin{equation*}
\| u_\delta\|_{H^1(B_R)} \le C_R \| f\|_{L^2},  
\end{equation*}
for some $C_R$ independent of $f$ and $\delta$. 
\end{proposition}

When $\Omega_j = B_{r_j}$ for $j=2,3$ and $A = I$ in $\Omega_3 \setminus \Omega_2$,   more quantitative estimates on the blow up and the  boundedness of the power are given in the following proposition in the spirit of \cite[Theorems 1.2 and 1.3]{Ng-CALR} (inspired by  \cite{AmmariCiraoloKangLeeMilton, KohnLu}). 

\begin{proposition} \label{pro3} Let $d=2, 3$, $f \in L^2_{\mc}(\mR^d)$, and  let $u_\delta \in H^1(\mR^d)$ be the unique solution of \eqref{def-ud}.  Assume that $(s_0A, s_0\Sigma)$ is doubly complementary media, $\Omega_2 = B_{r_2}$ and $\Omega_3 = B_{r_3}$ for  some $0 < r_2 < r_3$,  and $(A, \Sigma) = (I, 1)$ in $B_{r_3} \setminus B_{r_2}$. We have 
\begin{enumerate}
\item If  there exists $w \in H^1(B_{r_0} \setminus B_{r_2})$ for some $r_0 > \sqrt{r_2 r_3}$ with the properties 
\begin{equation*}
\dive (A \nabla w) + k^2 \Sigma w = f \mbox{ in } B_{r_0} \setminus B_{r_2} \quad \mbox{ and } \quad  \quad w = A \nabla w \cdot \nu = 0 \mbox{ on } \partial B_{r_2},
\end{equation*}
then 
\begin{equation*}
\limsup_{\delta \to 0} \delta \|u_\delta \|_{H^{1}(B_{r_3})}^2 < + \infty. 
\end{equation*}

\item If  there does {\bf not} exist  $v \in H^1(B_{r_0} \setminus B_{r_2})$ for some $r_0 <  \sqrt{r_2 r_3}$ with the properties 
\begin{equation*}
\dive (A \nabla v) + k^2 \Sigma v = f \mbox{ in } B_{r_0} \setminus B_{r_2} \quad \mbox{ and } \quad  \quad v = A \nabla v \cdot \nu = 0 \mbox{ on } \partial B_{r_2},
\end{equation*}
then 
\begin{equation*}
\liminf_{\delta \to 0} \delta \| \nabla u_\delta \|_{L^2(B_{r_3} \setminus B_{r_2})}^2= + \infty. 
\end{equation*}
\end{enumerate}
\end{proposition}

One only assumes that $(A, \Sigma) = (I, 1)$ in $B_{r_3} \setminus B_{r_2}$: the separation of variables is out of reach here. The  proof of the first statement of  Proposition~\ref{pro3} is based on a kind of removing singularity technique and has roots from \cite{Ng-CALR}. A key point is the construction of the auxiliary function $W_\delta$ in \eqref{def-Vd-S}. The proof of the second statement is  based on an observation  on a Cauchy problem in \cite{Ng-CALR} and involves a three spheres inequality. 

\medskip 
We finally point out that the stability of the Helmholtz with sign changing coefficients was studied in \cite{AnneSophieChesnelCiarlet1, CostabelErnst, KettunenLassas} and references therein by the integral method, the pseudo differential operator theory, and the T-coercivity approach and was unified and generalized in  \cite{Ng-WP} using different techniques based on reflections and the study of Cauchy's problems. It was also showed in \cite{Ng-WP} that complementary property is necessary for the appearance of resonance. 

\medskip 
The paper is organized as follows. The proof of Theorem~\ref{thm-main} is given in Section~\ref{sect-thm-main}. Sections~\ref{sect-pro1}, \ref{sect-pro2},  and \ref{sect-pro3} are devoted to the proofs of  Propositions~\ref{pro1}, \ref{pro2}, and \ref{pro3} respectively. Finally, in Section~\ref{sect-cloaking}, we present schemes of cloaking a general class of  sources via ALR in arbitrary medium for the finite frequency regime. 

\section{Proof of Theorem~\ref{thm-main}}\label{sect-thm-main}

We start this section with a lemma on the stability of \eqref{def-ud} which is used repeatedly in this paper.  

\begin{lemma}\label{lem1} Let $d=2, \, 3$, $R_0> 0$, $k>0$,  $\delta \in (0, 1)$, $R_0> 0$, $g \in H^{-1}(\mR^d)$  with support in $B_{R_0}$. Then there exists a unique outgoing  solution $v_\delta \in H^1_{\loc}(\mR^d)$ to the equation
\begin{equation}\label{eqd}
\dive (\epss_\delta A \nabla v_\delta) + k^2 \epss_0 \Sigma v_\delta = g  \mbox{ in } \mR^d. 
\end{equation}
Moreover,
\begin{equation}\label{stability}
\| v_\delta \|_{H^1(B_R)}^2 \le  \frac{C_R}{\delta} \Big| \int g \bar v_\delta \Big|, 
\end{equation}
for some positive constant $C_R$ independent of $g$ and $\delta$, as $\delta$ is small.
\end{lemma}

\noindent{\bf Proof.} The existence of $v_\delta$ follows from the uniqueness of $v_{\delta}$ by using the limiting absorption principle. We now establish the uniqueness of $v_{\delta}$ by showing that $v_\delta = 0$ if $v_\delta \in H^1_{\loc}(\mR^d)$ is an outgoing  solution  of 
\begin{equation*}
\dive(\epss_\delta A \nabla v_\delta) + k^2 \epss_0 \Sigma v_\delta = 0 \mbox{ in } \mR^d.
\end{equation*}
Multiplying the above equation by $\bar v_\delta$ (the conjugate of $v_{\delta}$) and integrating the obtained expression on $B_R$ with $R \ge R_0$, we have
\begin{equation*}
\int_{B_R} \epss_\delta \langle A \nabla v_\delta, \nabla v_\delta \rangle \, dx - \int_{B_R} k^2 \epss_0 \Sigma |v_\delta|^2 \, dx - \int_{\partial B_R} \partial_r v_\delta \bar v_\delta = 0.
\end{equation*}
By letting $R \to + \infty$, considering the imaginary part, and using the outgoing condition, we obtain 
\begin{equation*}
\int_{\Omega_2 \setminus \Omega_1} \langle A \nabla v_\delta, \nabla v_\delta \rangle \, dx= 0.
\end{equation*}
Since $A$ is uniformly elliptic, it follows that $v_\delta$ is constant in $\Omega_2 \setminus \Omega_1$. Thus $v_\delta = 0$ in $\Omega_2 \setminus \Omega_1$ since
$\dive(\epss_\delta A \nabla v_\delta) + k^2 \epss_0 \Sigma v_\delta = 0$ in $\Omega_2 \setminus \Omega_1$.
This implies $v_\delta = 0 $ in $\mR^d \setminus \Omega_{2}$ and $v_\delta = 0$ in $\Omega_{1}$ by the unique continuation principle. The proof for the uniqueness of $v_\delta$ is complete.

\medskip

Without loss of generality, one may assume that such that \eqref{cond-I} holds for $|x| \ge R_0$ and $\Omega_2 \subset \subset B_{R_0}$.  We next establish \eqref{stability} with $R = R_0$  by  contradiction. Assume that \eqref{stability} is not true. Then 
there exists $(g_{\delta}) \subset H^{-1}(\mR^d)$ such that
\begin{equation}\label{contradict-assumption}
\| v_\delta\|_{H^1(B_{R_0})} =1 \mbox{ and }\frac{1}{\delta} \Big| \int g_\delta \bar v_\delta \Big|   \to 0,
\end{equation}
as $\delta \to 0$, where $v_\delta \in H^1_{\loc}(\mR^d)$ is the unique solution to the equation
\begin{equation}\label{eq-v-delta}
\dive (\epss_\delta A \nabla v_\delta) + k^2 \epss_0 \Sigma v_\delta = g_\delta  \mbox{ in } \mR^d.
\end{equation}
In fact, by contradiction these properties only hold for a sequence of $(\delta_n) \to 0$. However, for simplicity of the notation, we still use $\delta$ instead of $\delta_n$  to denote an element of such a sequence. Since (see e.g., \cite[Lemma 2.3]{NguyenHelmholtz})
\begin{equation}\label{est-Helm}
\| v_\delta \|_{H^1(B_R \setminus B_{R_0})} \le C_R \| v_\delta\|_{H^{1/2}(\partial B_{R_0})}, 
\end{equation}
without loss of generality, one may assume that $v_\delta \to v_0$ in $L^2_{\loc}(\mR^d)$ and weakly in $H^1_{\loc}(\mR^d)$ for some $v_0 \in H^1_{\loc}(\mR^d)$. 
Multiplying  \eqref{eq-v-delta} by $\bar v_\delta$ and integrating the obtained expression in $B_R$ with $R\ge R_0$, we have
\begin{equation}\label{part1-lem1}
\int_{B_R} \epss_\delta \langle A \nabla v_\delta, \nabla v_\delta \rangle  \, dx -  \int_{B_R} k^2 \epss_0 \Sigma |v_\delta|^2 \, dx = - \int_{B_R} g_\delta \bar v_\delta \, dx + \int_{\partial B_{R}} \partial_r v_\delta \bar v_\delta.
\end{equation}
Letting $\delta \to 0$, by \eqref{contradict-assumption}, we obtain 
\begin{equation}\label{part2-lem1}
\Im \Big( \int_{\partial B_{R}} \partial_r v_0 \bar v_0 \Big) =0.
\end{equation}
Using the fact that $v_0$ satisfies the outgoing condition, by  letting $R \to + \infty$ and considering the imaginary part in \eqref{part2-lem1}, we obtain 
\begin{equation*}
\int_{\partial B_R} |v_0|^2 \to 0 \mbox{ as } \delta \to 0. 
\end{equation*}
It follows from Reillich's lemma that $v_0= 0$ in $\mR^d \setminus B_{R_0}$. One derives from  the unique continuation principle that 
\begin{equation}\label{toto1-v0}
v_0 = 0 \mbox{ in } \mR^d. 
\end{equation}
Letting $R \to \infty$, considering the imaginary part in \eqref{part1-lem1}, and using \eqref{contradict-assumption},  we obtain 
\begin{equation}\label{toto1-nabla}
\| \nabla v_\delta\|_{L^2(\Omega_2  \setminus \Omega_1)} \to 0 \mbox{ as } \delta \to 0.
\end{equation}
Considering the real part of \eqref{part1-lem1} with $R= R_0$, we derive from \eqref{toto1-v0} and \eqref{toto1-nabla} that
\begin{equation*}
\| v_{\delta}\|_{H^{1}(B_{R_0})} \to 0 \mbox{ as } \delta \to 0.  
\end{equation*}
We have a contradiction by \eqref{contradict-assumption}. Hence \eqref{stability} holds for $R = R_0$. The conclusion now follows from \eqref{est-Helm}. \proofend

\medskip 
We are ready to give 

\medskip 
\noindent{\bf Proof of Theorem~\ref{thm-main}.} Define 
\begin{equation*}
u_{1, \delta} = u_\delta \circ F^{-1} \mbox{ in } \mR^d \setminus \Omega_2, 
\end{equation*}
and 
\begin{equation*}
u_{2, \delta} = u_{1, \delta} \circ G^{-1} \mbox{ in } \Omega_3. 
\end{equation*}
It follows from  \eqref{def-DC} and a change of variables (see e.g., \cite[Lemma 2]{Ng-Complementary}) that 
\begin{equation}\label{eq-u12}
\dive (A \nabla u_{1, \delta})  + k^2 \Sigma u_{1, \delta} + i \delta \dive \big( A \nabla u_{1, \delta} \big) = \dive(A \nabla u_{2, \delta}) + k^2 \Sigma u_{2, \delta} = 0 \mbox{ in } \Omega_3 \setminus \Omega_2. 
\end{equation}
Set
\begin{equation}\label{def-Un-cor}
\hu_\delta = \left\{\begin{array}{cl} 
u_\delta & \mbox{ in } \mR^d \setminus \Omega_3, \\[6pt]
u_\delta - (u_{1, \delta} - u_{2, \delta}) & \mbox{ in } \Omega_3  \setminus \Omega_2, \\[6pt]
u_{2, \delta} & \mbox{ in } \Omega_2. 
\end{array}\right.
\end{equation}
It follows from \eqref{eq-u12} that $\hu_\delta \in H^1_{\loc}(\mR^d)$ is the unique outgoing solution of 
\begin{equation}\label{eq-Un-cor}
\dive (\hat A \nabla \hu_\delta) + k^2 \hat \Sigma \hu_{\delta} = f + \dive \big([s_\delta - s_0] A \nabla u_{1, \delta} \big)   \mbox{ in } \mR^d. \end{equation}
Since,  by Lemma~\ref{lem1}, 
\begin{equation}\label{est1-TM}
\| u_\delta\|_{H^1(B_{R})} \le C_R \delta^{-1} \|f \|_{L^2}. 
\end{equation}
It follows from \eqref{eq-Un-cor} and Lemma~\ref{lem1} that 
 \begin{equation*}
\|\hu_\delta \|_{H^1(B_{R})} \le C_R \| f\|_{L^2}. 
\end{equation*}
As a consequence, we have 
\begin{equation}\label{est4-TM}
\| u_\delta\|_{H^1(B_R \setminus \Omega_3)} \le C_R \| f\|_{L^2}.  
\end{equation}
Applying \cite[Theorem 5.3]{AlessandriniRondi} (a consequence of a three spheres inequality) for $u_\delta$ in $B_{R} \setminus \Omega_2$ for some $R>0$ such that $\Omega_3 \subset \subset B_R$, we derive from \eqref{est1-TM} and \eqref{est4-TM} that 
\begin{equation}\label{est5-TM}
\lim_{\delta \to 0} \delta \| u_\delta\|_{L^2(B_R)} = 0. 
\end{equation}
Using \eqref{stability} in Lemma~\ref{lem1}, we derive from \eqref{est5-TM} that
\begin{equation*}
\lim_{\delta \to 0} \delta \| u_\delta\|_{H^1(B_R)} = 0; 
\end{equation*} 
which implies
\begin{equation}\label{est2-1-TM}
\lim_{\delta \to 0} \| \dive \big([s_\delta - s_0] A \nabla u_{1, \delta} \big) \|_{H^{-1}} = 0. 
\end{equation}
A combination of \eqref{eq-Un-cor} and  \eqref{est2-1-TM} yields, for $R> 0$,  
\begin{equation*}
\|\hu_\delta - \hu \|_{H^1(B_{R})}  \to 0 \mbox{ as } \delta \to 0. 
\end{equation*}
The proof is complete. 
\proofend

\begin{remark} \fontfamily{m} \selectfont The same proof gives the same conclusion even if one replaces $f$ by  $f_\delta$ in \eqref{def-ud} and assumes that $(f_\delta)$ converges to $f$ in $L^2(\mR^d)$ and $\supp f_\delta \subset \subset B_{R_0} \setminus \Omega_2$ for some $R_0 > 0$.  
\end{remark}

\begin{remark} \fontfamily{m} \selectfont One of the key points in the proof is the definition of $\hu_\delta$  in \eqref{def-Un-cor} after introducing $u_{1, \delta}$  and $u_{2, \delta}$ as in \cite{Ng-CALR}. In $\Omega_3 \setminus \Omega_2$, we remove $u_{1, \delta} - u_{2, \delta}$ from $u_\delta$. The removing term is the singular part of $u_\delta$ in $\Omega_3 \setminus \Omega_2$. 
The way of defining the removing term is intrinsic and more robust than the one in \cite{Ng-CALR}, which is based on the separation of variables for a general shell developed there. As  seen from there,  the removing term becomes more and more singular when one approaches $\partial \Omega_2$. The idea of removing the singular term was inspired from the study of  the Ginzburg-Landau equation in the work of Bethuel, Brezis, and Helein in \cite{BBH}.  Another new important point in the proof is to establish \eqref{est5-TM}. This is obtained by first proving  that $u_\delta$ is bounded outside $\Omega_3$ (again based on the behaviour of $\hat u_\delta$) and  then applying a three spheres inequality. 
\end{remark}

\section{Proof of Proposition~\ref{pro1}} \label{sect-pro1}

We prove Proposition~\ref{pro1} by contradiction. Assume that 
\begin{equation}\label{contradiction-pro1}
\limsup_{\delta \to 0} \delta \| \nabla u_\delta\|_{L^2(\Omega_2 \setminus \Omega_1)}^2 < + \infty. 
\end{equation}
Since $\dive(A \nabla u_\delta) + k^2 s_0 s_\delta^{-1} \Sigma u_\delta  = 0$ in $\Omega_2 \setminus \Omega_1$, it follows from a compactness argument that 
\begin{equation*}
\| u_\delta\|_{L^2(\Omega_2 \setminus \Omega_1)} \le C \|\nabla u_\delta \|_{L^2(\Omega_2 \setminus \Omega_1)}. 
\end{equation*}
We derive that 
\begin{equation}\label{contradiction1-pro1}
\lim_{\delta \to 0} \delta \|u_\delta \|_{H^1(\Omega_2 \setminus \Omega_1)}^2 < + \infty; 
\end{equation}
this implies, for $R> 0$,  
\begin{equation*}
\lim_{\delta \to 0} \delta \|u_\delta \|_{H^1(B_R)}^2 = 0. 
\end{equation*}
Define 
\begin{equation*}
u_{1, \delta} = u_{\delta} \circ F^{-1} \mbox{ in } \Omega_3 \setminus \Omega_2 
\end{equation*}
and set 
\begin{equation*}
v_\delta = u_{1, \delta} - u_\delta \mbox{ in } D. 
\end{equation*}
Then 
\begin{equation}\label{part1-pro1}
\dive(A \nabla v_\delta ) + k^2 \Sigma v_{\delta} = g_\delta \mbox{ in } D, 
\end{equation}
\begin{equation}\label{part2-pro1}
v_\delta = 0  \mbox{ on } D \cap \partial \Omega_2 \quad \mbox{ and } \quad A \nabla v_\delta \cdot \nu = h_\delta \mbox{ on } \partial D \cap \partial \Omega.  
\end{equation}
Here 
\begin{equation*}
g_\delta = - i \delta \dive \big(A \nabla u_{1, \delta} \big) =\frac{i \delta}{1 + i \delta} k^2 \Sigma u_{1, \delta} \mbox{ in } D
\end{equation*}
and 
\begin{equation*}
h_\delta = i \delta \nabla u_{1, \delta} \cdot \nu \mbox{ on } \partial D \cap \partial \Omega_2. 
\end{equation*}
It is clear from \eqref{contradiction-pro1} that 
\begin{equation}\label{part3-pro1}
\delta^{1/2}\|g_\delta \|_{L^2(D)} +  \delta^{1/2}\| h_\delta\|_{H^{-1/2}(\partial D \cap \partial \Omega_2)} \le C,  
\end{equation}
for some $C$ independent of $\delta$. 
Using \eqref{part1-pro1}, \eqref{part2-pro1}, and \eqref{part3-pro1}, and applying \cite[Lemma 10]{Ng-WP}, we have 
\begin{equation*}
\limsup_{\delta \to 0} \delta^{1/2}\| v_\delta\|_{H^1(D)} = + \infty:  
\end{equation*}
which contradicts \eqref{contradiction1-pro1}. The proof is complete. \proofend

\section{Proof of Proposition~\ref{pro2}} \label{sect-pro2}

We in fact prove a slightly more general result. 

\begin{proposition}\label{pro-compatible} Let $d=2, 3$, $\delta \in (0, 1)$, $f \in L^2_{\mc}(\mR^d)$, $g \in H^{1/2}(\partial \Omega_3)$,  and $h \in H^{-1/2}(\partial \Omega_3)$. Assume that  $(s_0 A, s_0 \Sigma)$ is doubly complementary and $\supp f \cap   \Omega_3 = \O$,  and  let $V_\delta \in H^1_{\loc}(\mR^d \setminus \partial \Omega_3)$ be the unique outgoing  solution of
\begin{equation}\label{eq-Vd}
\left\{\begin{array}{c}\dive (s_\delta A \nabla V_\delta)  + k^2 s_0 \Sigma V_\delta = f \mbox{ in } \mR^d \setminus \partial \Omega_3,  \\[6pt]
[V_\delta] = g \quad \mbox{ and } \quad [A \nabla V_\delta \cdot \eta] = h \mbox{ on } \partial \Omega_3. 
\end{array}\right.
\end{equation}
Then, for $R > 0$,  
\begin{equation*}
 \|V_\delta \|_{H^{1}(B_R \setminus \partial \Omega_3)} \le C_R \big(\| f\|_{L^2(\Omega)} + \| g\|_{H^{1/2}(\partial \Omega_3)} + \| h\|_{H^{-1/2}(\partial \Omega_3)}\big),  
\end{equation*}
for some positive constant $C_R$ independent of $\delta$, $f$, $g$, and $h$. 
\end{proposition}

Here and in what follows in this paper, we denote $[v]  = v \Big|_{ext} - v \Big|_{int}$  and $[M \nabla v \cdot \eta] = M \nabla v \cdot \nu  \Big|_{ext} - M \nabla v \cdot \nu \Big|_{ext} $ on $\partial \Omega$ for a smooth bounded open subset $\Omega$ of $\mR^d$, a matrix-valued function $M$, and for an appropriate function $v$.

\medskip 
It is clear that Proposiiton~\ref{pro-compatible} implies Proposition~\ref{pro2}. Proposition~\ref{pro-compatible}  is also used in the proof of Proposition~\ref{pro3}.

\medskip 

\noindent{\bf Proof.} The proof of Proposition~\ref{pro-compatible} has roots from \cite{Ng-Complementary} and the key point is to construct a solution $V_0$ of \eqref{eq-Vd} for $\delta = 0$. Let $\cV \in H^1_{\loc}(\mR^d \setminus \partial \Omega_3)$ be the unique outgoing solution to 
\begin{equation*}
\left\{\begin{array}{c}\dive (\cA \nabla \cV) + k^2 \cS \cV = f \mbox{ in } \mR^d \setminus \partial \Omega_3,  \\[6pt]
[\cV] = g \quad \mbox{ and } \quad [\cA \nabla \cV \cdot \nu] = h \mbox{ on } \partial \Omega_3,  
\end{array}\right.
\end{equation*}
where $(\cA, \cS)$ is defined in \eqref{def-cAS}. 
Then
\begin{equation}\label{haha1-Thm2.3}
\| \cV \|_{H^1(B_{R} \setminus \partial \Omega_3)} \le C_R \big(\| f\|_{L^2} + \| g\|_{H^{1/2}(\partial \Omega_3) } + \| h\|_{H^{-1/2}(\partial \Omega_3)}\big).  
\end{equation}
Define $V_0 \in H^1(\Omega \setminus \partial B_{r_3})$ as follows 
\begin{equation}\label{haha2-Thm2.3}
V_0 =  \left\{\begin{array}{cl} \cV &\mbox{ in } \mR^d \setminus \Omega_2,  \\[6pt]
\cV \circ F & \mbox{ in } \Omega_2 \setminus \Omega_1, \\[6pt]
\cV \circ G \circ F & \mbox{ in } \Omega_1. 
\end{array}\right.
\end{equation}
Using \eqref{def-DC} and applying \cite[Lemma 2]{Ng-Complementary} (a change of variables), as in \cite[Step 2 in Section 3.2.2]{Ng-Complementary}, one can verify that $V_0 \in H^1_{\loc}(\mR^d \setminus \partial \Omega_3)$ is an outgoing solution to 
\begin{equation*}
\dive (s_0 A \nabla V_0) + k^2 s_0 \Sigma V_0= f \mbox{ in } \mR^d \setminus (\partial \Omega_3 \cup \partial \Omega_1). 
\end{equation*}
One also obtains from the definition of $V_0$ and $\hat V$ that  
\begin{equation*}
[V_0] = g \quad \mbox{ and } \quad [A \nabla V_0 \cdot \nu] = h \mbox{ on } \partial \Omega_3. 
\end{equation*}
and 
\begin{equation*}
[V_0] = 0 \quad \mbox{ and } \quad [A \nabla V_0 \cdot \nu] = 0 \mbox{ on } \partial \Omega_1. 
\end{equation*}
Hence $V_0 \in H^1_{\loc}(\mR^d)$ is an outgoing solution of \eqref{eq-Vd} with $\delta = 0$.  Set  
\begin{equation}\label{def-Wd-Thm2.3}
W_\delta = V_\delta - V_0 \mbox{ in } \Omega. 
\end{equation}
Then $W_\delta \in H^1_0(\Omega)$ is the unique solution to 
\begin{equation*}
\dive(s_\delta A \nabla W_\delta)  + k^2 s_0 \Sigma W_\delta =  - \dive\big( i \delta A \nabla V_0  1_{B_{r_2} \setminus B_{r_1} } \big)  \mbox{ in } \mR^d.
\end{equation*}
Here and in what follows,  for a subset $D$ of $\mR^d$, $1_D$ denotes the characteristic function of $D$. 
Applying Lemma~\ref{lem1}, we have
\begin{equation}\label{haha3-Thm2.3}
\|W_\delta\|_{H^1(B_{R})} \le C_R \|V_0\|_{H^1(\Omega_2 \setminus \Omega_1)}. 
\end{equation}
The conclusion follows from \eqref{haha1-Thm2.3}, \eqref{haha2-Thm2.3}, \eqref{def-Wd-Thm2.3}, and \eqref{haha3-Thm2.3}. \proofend

\section{Proof of Proposition~\ref{pro3}}\label{sect-pro3}

\noindent{\bf Step 1:} Proof of the first statement. Without loss of generality, one might assume that $r_2 = 1$. Define 
\begin{equation*}
u_{1, \delta} = u_\delta \circ F^{-1} \mbox{ in } \mR^d \setminus B_{r_3}, 
\end{equation*}
and
\begin{equation*}
u_{2, \delta} = u_{1, \delta} \circ G^{-1} \mbox{ in } B_{r_3}. 
\end{equation*}
Let $\phi \in H^1(B_{r_3} \setminus B_{r_2})$ be the unique solution to 
\begin{equation}\label{def-phi-Thm2.3}
\Delta \phi + k^2 \phi = f \mbox{ in } B_{r_3} \setminus B_{r_2},  \quad \phi = 0 \mbox{ on } \partial B_{r_2}, \quad \mbox{ and } \quad \partial_r \phi - i k \phi = 0 \mbox{ on } \partial B_{r_3}, 
\end{equation}
and set 
\begin{equation*}
W = w - \phi \mbox{ in } B_{r_0} \setminus B_{r_2}. 
\end{equation*}
Then $W \in H^1(B_{r_3} \setminus B_{r_2})$ satisfies 
\begin{equation}\label{pro-V-S} 
\Delta W + k^2 W = 0 \mbox{ in } B_{r_0} \setminus B_{r_2}, \quad W = 0 \mbox{ on } \partial B_{r_2}, \quad \mbox{ and } \quad \partial_r W = - \partial_r \phi \mbox{ on } \partial B_{r_2}. 
\end{equation}

We now consider the case $d=2$ and $d=3$ separately. 

\medskip 

\noindent \underline{Case 1:} $d=2$.  As in \cite{Ng-Superlensing}, define
\begin{equation}\label{def-Jn}
\hat J_n(r) = 2^n n! J_n(r) \quad \mbox{ and } \quad \hat Y_n(r) = \frac{\pi i}{2^{n} (n-1)!} Y_n(r), 
\end{equation}
where $J_n$ and $Y_n$ are the Bessel and Neumann functions of order $n$. It follows from \cite[(3.57) and (3.58)]{ColtonKressInverse} that 
\begin{equation}\label{bh1}
\hat J_n (t)  = t^{n}\big[1 + o(1) \big]
\end{equation}
and
\begin{equation}\label{bh2}
\hat Y_n (t)  = t^{-n} \big[1 + o(1) \big], 
\end{equation}
as $n \to + \infty$.  Since $\Delta W + k^2 W = 0$ in $B_{r_3} \setminus B_{r_2}$, 
one can represent $W$ as follows
\begin{equation}\label{re-u-k}
W =     \sum_{n=0}^\infty \sum_{\pm} \big[a_{n, \pm} \hat J_{n}(|x|) + b_{n, \pm} \hat Y_{n}(|x|) \big] e^{\pm i n \theta} \quad \mbox{ in } B_{r_3} \setminus B_{r_2},
\end{equation}
for $a_{n, \pm}, b_{n, \pm} \in \mC$ ($n \ge 0$) with $a_{0, + } = a_{0, -}$ and $b_{0, +} = b_{0, -}$.  Since $r_2 = 1 <  r_0$ and $W=0$ on $\partial B_{r_2}$,   we derive from \eqref{bh1} and \eqref{bh2} that, for some $N>0$ independent of $W$, 
\begin{equation}\label{finite-V}
\| W\|_{H^1(B_{r_0} \setminus B_{r_2})}^2 \sim \sum_{n = 0}^N \sum_{\pm} (|a_{n, \pm}|^2 + |b_{n, \pm}|^2 ) + \sum_{n = N+1}^\infty \sum_{\pm} n |a_{n, \pm} ^2| r_0^{2 n} < + \infty. 
\end{equation}
Here we used the fact $W=0$ on $\partial B_{r_2}$ to derive that $a_{n, \pm} \sim b_{n, \pm}$ for $n \ge N+1$ since $r_2 = 1$. 
One of the keys  in the proof is the construction of  $W_\delta \in H^1(B_{r_3} \setminus B_{r_2})$ which is  defined as follows 
\begin{equation}\label{def-Vd-S}
W_\delta  =   \sum_{n=0}^\infty \sum_{\pm} \frac{1}{1 + \xi_{n}}\big[a_{n, \pm} \hat J_{n}(|x|) + b_{n, \pm} \hat Y_{n}(|x|) \big] e^{\pm i n \theta}  \mbox{ in } B_{r_3} \setminus B_{r_2}, 
\end{equation}
where 
\begin{equation}\label{xil}
\xi_n =  \delta^{1/2} (r_3 / r_0)^{n} \mbox{ for } n \ge 0. 
\end{equation}
Roughly speaking, $W_\delta$ is   the singularity of $u_\delta$.
From the definition of $W_\delta$, we have
\begin{equation}\label{eq-Wdd}
\Delta W_\delta  + k^2 W_{\delta}= 0 \mbox{ in } B_{r_3} \setminus \bar B_{r_2},  \quad W_\delta = 0 \mbox{ on } \partial B_{r_2}, 
\end{equation}
and 
\begin{equation}\label{vdelta-H1}
\|W_\delta \|_{H^1(B_{r_3} \setminus B_{r_2})}^2 \sim\sum_{n = 0}^N \sum_{\pm} (|a_{n, \pm}|^2 + |b_{n, \pm}|^2 )   +  \sum_{n= N+1 }^\infty  \sum_{\pm} \frac{n |a_{n, \pm}|^2}{1 + \xi_n^2} r_3^{2 n}. 
\end{equation}
By \eqref{xil}, we have, if $\xi_n \le 1$, then
\begin{equation}\label{control1-v}
\frac{n |a_{n, \pm}|^2}{1 + \xi_n^2}  r_3^{2 n} \le 
 n   |a_{n, \pm}|^2r_3^{2n} \le \delta^{-1} n  |a_{n, \pm}|^2 r_0^{2n}, 
\end{equation}
and if $\xi_n \ge 1$, then
\begin{equation}\label{control2-v}
\frac{ n|a_{n, \pm}|^2  }{1 + \xi_n^2} r_3^{2 n} \le n |a_{n, \pm}|^2 r_3^{2 n} \xi_n^{-2} 
= \delta^{-1} n |a_{n, \pm}|^2 r_0^{2n}. 
\end{equation}
A combination of \eqref{finite-V},  \eqref{vdelta-H1}, \eqref{control1-v}, and \eqref{control2-v} yields  
\begin{equation}\label{vdelta-H1-1}
\|W_\delta \|_{H^1(B_{r_3} \setminus B_{r_2})} \le C \delta^{-1/2}. 
\end{equation}
Let $W_{1, \delta} \in H^1_{\loc}(\mR^2)$ be the unique outgoing solution to 
\begin{equation*}\left\{
\begin{array}{c}\dive(s_\delta A \nabla W_{1, \delta}) + k^2 s_0 \Sigma W_{1, \delta}= 0 \mbox{ in } \mR^2 \setminus \partial B_{r_2}, \\[6pt]
 [s_\delta A \nabla W_{1, \delta}\cdot \nu] = (-1 - i \delta) h_\delta  \mbox{ on } \partial B_{r_2},
\end{array}\right. 
\end{equation*}
where 
\begin{equation*}
h_\delta = - \partial_r (\phi + W_\delta)  \mbox{ on } \partial B_{r_2}, 
\end{equation*}
and let $W_{2, \delta} \in H^1_{\loc}(\mR^2  \setminus \partial B_{r_3})$ be the unique outgoing solution to 
\begin{equation*}\left\{
\begin{array}{c}
\dive(s_\delta A \nabla W_{2, \delta})  +  k^2 s_0 \Sigma W_{2, \delta} = f 1_{\mR^2 \setminus B_{r_3}} \mbox{ in } \mR^2 \setminus \partial B_{r_3}, \\[6pt]
[W_{2,\delta}] = \phi +  W_\delta \quad \mbox{ and } \quad  [A \nabla W_{2, \delta} \cdot \nu] = \partial_r \phi  +  \partial_r W_\delta   \mbox{ on } \partial B_{r_3}.  
\end{array}\right.
\end{equation*}
Recall that,  for a subset $D$ of $\mR^d$, $1_D$ denotes the characteristic function of $D$. 
From \eqref{def-phi-Thm2.3}, \eqref{eq-Wdd}, and the fact $(A, \Sigma)=(I, 1)$ in $B_{r_3} \setminus B_{r_2}$, we have 
\begin{equation}\label{u-V-delta}
u_\delta - ( \phi + W_\delta) {\bf 1}_{B_{r_3} \setminus B_{r_2}} = W_{1, \delta} + W_{2, \delta} \mbox{ in } \Omega. 
\end{equation}
Using \eqref{pro-V-S} and \eqref{def-Vd-S}, we obtain, on $\partial B_{r_2}$, 
\begin{equation*}
h_\delta = - \partial_r (\phi + W_\delta) = \partial_r (W- W_\delta)  = \partial_r \left( \sum_{n= 0}^{\infty} \sum_{\pm} \frac{\xi_n}{1 + \xi_n}  \big[a_{n, \pm} \hat J_{n}(|x|) + b_{n, \pm} \hat Y_{n}(|x|) \big]e^{\pm i n \theta} \right). 
\end{equation*}
Since $r_2 = 1$, it follows that 
\begin{equation}\label{hdelta-H12-S}
\| h_\delta \|_{H^{-1/2}(\partial B_{r_2})}^2 \sle \sum_{n = 0}^N \sum_{\pm} (|a_{n, \pm}|^2 + |b_{n, \pm}|^2 ) + \sum_{n = N+1}^\infty \sum_{\pm} \frac{n |\xi_n|^2 }{1  + |\xi_n|^2} |a_{n, \pm}| ^2. 
\end{equation}
Using \eqref{xil}, we have, if $\xi_n \le 1$ then 
\begin{equation}\label{control1-h-S}
\frac{n |\xi_n|^2}{1 + |\xi_n|^2}   \le 
\delta n  |a_{n, \pm}|^2 (r_3/ r_0)^{2 n} = \delta  n r_0^{2n} (r_3/ r_0^2)^{2 n} \le \delta n r_0^{2n},  
\end{equation}
since $r_0 >  \sqrt{r_2 r_3} = \sqrt{r_3}$,  and if $\xi_n \ge 1$  then 
\begin{equation}\label{control2-h-S}
\frac{ n |\xi_n|^2 }{1 + |\xi_n|^2} \le n |a_{n, \pm}|^2  = n  r_0^{2 n } r_0^{-2 n} \le  \delta n  r_0^{2 n }, 
\end{equation}
since $\delta^{1/2} r_0^{n} > \delta^{1/2} (r_3/ r_0)^{n} \ge 1$. 
A combination of  \eqref{hdelta-H12-S}, \eqref{control1-h-S}, and \eqref{control2-h-S} yields  
\begin{equation*}
\| h_\delta \|_{H^{-1/2}(\partial B_{r_2})} \le C \delta^{1/2} \| W\|_{H^{1/2}(\partial B_{r_0})} \le C \delta^{1/2}.
\end{equation*}
Applying Lemma~\ref{lem1}, we have
\begin{equation}\label{W1-1}
\|W_{1, \delta} \|_{H^1(\Omega)} \le (C/ \delta) \delta^{1/2} = C\delta^{-1/2}. 
\end{equation}
On the other hand,  from \eqref{vdelta-H1-1} and Proposition~\ref{pro-compatible}, we obtain
\begin{equation}\label{W2-1}
\|W_{2, \delta} \|_{H^1(B_{r_3} \setminus B_{r_3})} \le C \delta^{-1/2}. 
\end{equation}
The conclusion in the case $d=2$ now follows from \eqref{vdelta-H1-1}, \eqref{u-V-delta}, \eqref{W1-1}, and \eqref{W2-1}. 

\medskip
\noindent \underline{Case 2:} $d=3$.   Define 
\begin{equation}\label{def-jn}
\hat j_n(t) =1 \cdot 3 \cdots (2n + 1) j_n(t) \quad \mbox{ and } \quad  \hat y_n = -  \frac{y_n(t)}{1 \cdot 3 \cdots (2n-1)} ,  
\end{equation}
where $j_n$ and $y_n$ are the spherical Bessel and Neumann functions of order n. 
Then, as $n$ large enough,  (see,  e.g., \cite[(2.37) and (2.38)]{ColtonKressInverse})  
\begin{equation}\label{jy-n}
\hat j_n(kr) = r^n \big(1 + O(1/n) \big) \quad \mbox{ and } \quad \hat y_n(kr) = r^{-n-1} \big(1 + O(1/n) \big). 
\end{equation}
Thus one can represent $W$ of the form
\begin{equation}\label{re-u-k1}
W=  \sum_{n=1}^\infty \sum_{-n }^n \big[a^n_m \hat j_{n}(|x|) + b^n_m \hat y_{n}(|x|) \big] Y^n_m(\hat x ) \quad \mbox{ in } B_{r_3} \setminus B_{r_0},
\end{equation}
for $ a^n_m, b^n_m \in \mC$ and $\hat x = x/ |x|$. Here $Y^n_m$ is the spherical harmonic function of degree $n$ and of order $m$.  Define $W_\delta \in H^1(B_{r_3} \setminus B_{r_2})$ by 
\begin{equation*}
W_\delta  =   \sum_{n=1}^\infty \sum_{-n }^n \frac{1}{1 + \xi_n} \big[a^n_m \hat j_{n}(|x|) + b^n_m \hat y_{n}(|x|) \big] Y^n_m(\hat x ) \mbox{ in } B_{r_3} \setminus B_{r_2}, 
\end{equation*}
where 
\begin{equation*}
\xi_n =  \delta^{1/2} (r_3 / r_0)^{n} \mbox{ for } n \ge 1. 
\end{equation*}
 The proof now follows similarly as in the case $d=2$. The details are left to the reader. 
 
 \medskip 
 \noindent{\bf Step 2}: Proof of the second statement.  Define $u_{1, \delta}  = u_{\delta} \circ F$ and  denote $u_{2^{-n}}$ and $u_{1, 2^{-n}}$ by $u_n$ and $u_{1, n}$. for notational ease.  We prove by contradiction that
\begin{equation}\label{claim-IE}
\limsup_{n \to +\infty} 2^{-n/2} \big( \| u_n\|_{H^{1}(B_{r_3} \setminus B_{r_2})}   + \| u_{1, n} \|_{H^{1}(B_{r_3} \setminus B_{r_2})}  \big)= + \infty. 
\end{equation} 
Assume that 
\begin{equation}\label{claim-IE-contradiction}
m: = \sup_{n}  2^{-n/2} \big( \| u_n\|_{H^{1}(B_{r_3} \setminus B_{r_2})}   + \| u_{1, n} \|_{H^{1}(B_{r_3} \setminus B_{r_2})}  \big)  <  + \infty. 
\end{equation} 
Define 
\begin{equation*}
v_n = u_n- u_{1, n} \mbox{ in } B_{r_3} \setminus B_{r_2} \quad \mbox{ and } \quad \phi_n =   i 2^{-n}  \partial_r  u_{1, n}  \mbox{ on } \partial B_{R_1}. 
\end{equation*}
Then 
\begin{equation*}
\Delta v_n + k^2 v_n  = f  \mbox{ in } B_{r_3} \setminus B_{r_2}, \; \;  v_n = 0 \mbox{ on } \partial B_{r_2}, \mbox{  and  } \quad \partial_r v_n =  \phi_n \mbox{ on } \partial B_{r_2}. 
\end{equation*}
We claim that $(v_n)$ is a Cauchy sequence in $H^1(B_{r_0} \setminus B_{r_2})$.  Indeed, set 
\begin{equation*}
V_{n} = v_{n+1} - v_{n}  \mbox{ in } B_{r_3} \setminus B_{r_2} \quad \mbox{ and } \quad \Phi_{n} = \phi_{n+ 1} - \phi_{n} \mbox{ on } \partial B_{r_2}. 
\end{equation*}
We have
\begin{equation*}
\Delta V_n + k^2 V_n= 0 \mbox{ in } B_{r_3} \setminus B_{r_2}, \quad V_{n}= 0 \mbox{ on } \partial B_{r_2}, \quad \mbox{ and } \quad \partial_r V_{n} = \Phi_{n} \mbox{ on } \partial B_{r_2}. 
\end{equation*}
From \eqref{claim-IE-contradiction}, we derive that 
\begin{equation*}
\| V_n\|_{H^1(B_{R_2} \setminus B_{R_1})} \le C m 2^{n/2} \quad \mbox{ and } \quad \| \Phi_n\|_{H^{1/2}(\partial B_{R_1})} \le C m 2^{-n/2}. 
\end{equation*}
In this proof, $C$ denotes a constant independent of $n$. Let $U_n \in H^1(B_{r_3})$ be the unique solution of 
\begin{equation*}
\Delta U_n + k^2 U_n = 0 \mbox{ in } B_{r_3} \setminus \partial B_{r_2},  \quad [\partial_r U_n] = \Phi_n, \quad \mbox{ and } \quad \partial_r U_n - i k U_n = 0 \mbox{ on } \partial B_{r_3}. 
\end{equation*}
We have 
\begin{equation}\label{pro-Un}
\| U_n\|_{H^1(B_{r_3})} \le C \|\Phi_n \|_{H^{-1/2}(\partial B_{r_2})}. 
\end{equation}
Applying Lemma~\ref{lem-threesphree} below (a three spheres inequality) for $V_n 1_{B_{r_3} \setminus B_{r_2}}- U_n$ in $B_{r_3}$,  we obtain from \eqref{pro-Un} that 
\begin{equation*}
\| V_{n} \|_{H^1(B_{r_0} \setminus B_{r_2})} \le C \Big( \| \Phi_n \|_{H^{-1/2}(\partial B_{r_2}) }^{\alpha} \| V_{n}\|_{H^1(B_{r_0} \setminus B_{r_2})}^{1 - \alpha} + \| \Phi_n \|_{H^{-1/2}(\partial B_{r_2}) } \Big) \le C m 2^{-n \beta},  
\end{equation*}
where $\alpha = \ln (r_3/ r_0) / \ln(r_3/ r_1) > 1/2$ and $\beta = \big(2 \alpha - 1 \big)/ 2 > 0$.  
Hence $(v_n)$ is a Cauchy sequence  in $H^1(B_{r_0} \setminus B_{r_2})$. 
Let $v $ be the  limit of $v_n$ in $H^1(B_{r_0} \setminus B_{r_2})$.  Then  
\begin{equation*}
\Delta v + k^2 v = f \mbox{ in } B_{r_0} \setminus B_{r_2}, \quad  v = 0 \mbox{ on } \partial B_{R_1}, \quad  \partial_r v  =  0 \mbox{ on } \partial B_{r}. 
\end{equation*}
This contradicts the non-existence of $v$.  Hence \eqref{claim-IE} holds. The proof is complete. \proofend

\medskip 
The following lemma is used in the proof of the second statement of Proposition~\ref{pro3}. 
\begin{lemma} \label{lem-threesphree} Let $d=2, 3$, $k, \, R>0$,   and let $u \in H^1(B_R)$ be a solution to the equation $\Delta u + k^2 u = 0$. Then, for $0 < R_1 < R_2 < R_3 \le R$, 
\begin{equation*}
\| u \|_{H^1(B_{R_2})} \le C_{R, k} \| u\|_{H^1(B_{R_1})}^{\alpha} \| u \|_{H^1(B_{R_3})}^{1 - \alpha}, 
\end{equation*}
where $\alpha =  \ln (R_3/ R_2) / \ln(R_3/ R_1)$ and $C_{R, k}$ is a positive constant independent of $R_1, R_2, R_3$, and $v$. 
\end{lemma}

\noindent{\bf Proof.} We first give the proof in two dimensions. Since $\Delta v + k^2 v = 0$ in $B_{R}$, one can represent $v$ of the form
\begin{equation*}
v=  \sum_{n=0}^\infty \sum_{\pm} a_{n, \pm} \hat J_{n}(|x|) e^{\pm i n \theta}  \mbox{ in } B_R, 
\end{equation*} 
for $a_{n, \pm} \in \mC$ ($n \ge 0$) with $a_{0, + } = a_{0, -}$ where $\hat J_n$ is defined in \eqref{def-Jn}. The conclusion now follows from \eqref{bh1} after applying H\"older's inequality.

The proof in three dimensions follows similarly. In this case, $v$ can be represented in the form 
\begin{equation*}
v=  \sum_{n=1}^\infty \sum_{-n }^n a^n_m \hat j_{n}(|x|) Y^n_m(\hat x ) \quad \mbox{ in } B_{R}, 
\end{equation*}
for $ a^n_m \in \mC$ and $\hat x = x/ |x|$ where $\hat j_n$ is defined in \eqref{def-jn}. The conclusion is now a consequence of \eqref{jy-n} after applying H\"older's inequality. \proofend

\section{Cloaking a source via anomalous localized resonance in the finite frequency regime}\label{sect-cloaking}

In this section, we describe how to use the theory discussed previously to cloak a source $f$ concentrating on an arbitrary bounded smooth manifold of codimension 1 in an arbitrary medium. We follow the strategy in \cite{Ng-CALR}.   
Without loss of generality, one may assume that the medium is contained in $B_{r_3} \setminus B_{r_2}$ for some $0< r_2 <  r_3$ and characterized by a matrix-valued function $a$ 
and a real bounded function $\sigma$. We assume in addition  that $a$ is Lipschitz and uniformly elliptic in $\overline{B_{r_3} \setminus B_{r_2}}$  and $\sigma$ is bounded below by a positive constant.  Let $f \in L^2(\Omega)$  for some bounded smooth open subset $\Omega \subset \subset B_{r_3} \setminus B_{r_2}$. One might assume as well that $\Omega \subset \subset B(x_0, r_0)$ for some $r_0 > 0$ and $x_0 \in \partial B_{r_2}$ where $r_0$ is the constant coming from Proposition~\ref{pro1}. 
Define $r_1 = r_2^2/ r_3$.  Let $F: B_{r_2} \setminus \{0 \} \to \mR^d \setminus B_{r_2}$ and $G: \mR^d \setminus B_{r_3} \to B_{r_3} \setminus \{0 \}$ be the Kelvin transform with respect to $\partial B_{r_2}$ and $\partial B_{r_3}$ respectively. Define 
\begin{equation}\label{A-cloak}
A, \Sigma= \left\{\begin{array}{cl}  a, \sigma & \mbox{ in } B_{r_3} \setminus B_{r_2}, \\[6pt]
F^{-1}_*a , F^{-1}_*\sigma  & \mbox{ in } B_{r_2} \setminus B_{r_1}, \\[6pt]
F^{-1}_* G^{-1}_*a,  F^{-1}_* G^{-1}_*\sigma & \mbox{ in } B_{r_1} \setminus B_{r_1^2/ r_2}, \\[6pt]
I, 1 & \mbox{ otherwise}.  
\end{array}\right.
\end{equation}
It is clear that $(s_0A, s_0 \Sigma)$ is doubly complementary. 
Applying Theorem~\ref{thm-main} and Proposition~\ref{pro1}, we obtain
\begin{proposition} \label{pro-cloaking} Let $d \ge 2$, $\delta > 0$,  and $\Omega \subset \subset  D: = B(x_0, r_0) \cap (B_{r_3} \setminus B_{r_2})$  be smooth and open, let 
$f \in L^2(\partial \Omega)$
and  let $u_\delta$ and $v_\delta$ be defined by \eqref{def-ud} and \eqref{def-vd} where $(A, \Sigma)$ is given in \eqref{A-cloak}. 
Assume that $f \not \in {\cal H}$ where 
\begin{equation*} 
{\cal H}: = \{A \nabla v \cdot \nu \big|_{\partial \Omega};  v \in H^1_0(\Omega) \mbox{ is a solution of } \dive(A \nabla v) + k^2 \Sigma v = 0 \mbox{ in } \Omega\}. 
\end{equation*}
There exists a sequence $\delta_n \to 0$ such that 
\begin{equation*}
\lim_{n \to \infty} E_{\delta_n} (u_{\delta_n}) = + \infty. 
\end{equation*}
Moreover, 
\begin{equation*}
v_{\delta_n} \to 0 \mbox{ weakly in } H^1_{\loc}(\mR^d \setminus B_{r_3}). 
\end{equation*}
\end{proposition} 

\noindent {\bf Proof.} By Theorem~\ref{thm-main} and Proposition~\ref{pro1}, it suffices to prove that there is no $W \in H^1(D)$ such that 
\begin{equation*}
\dive(A \nabla W) + k^2 \Sigma = f \mbox{ in } D  \quad \mbox{ and } \quad W= A \nabla W \cdot \eta = 0 \mbox{ on } \partial D \cap \partial B_{r_2}. 
\end{equation*}
In fact, Theorem~\ref{thm-main} and Proposition~\ref{pro1} only deal with the case $f \in L^2$, however, the same results hold for $f$  stated here and  the proofs are unchanged. Suppose that this is not true, i.e., such a $W$ exists.
Since $\dive(A \nabla W )  + k^2 \Sigma W = 0 $ in $D \setminus \bar \Omega$ and $W = A \nabla W \cdot \nu = 0$ on $\partial D \cap \partial B_{r_2}$,  it follows from the unique continuation principle that 
$W = 0$ in $D \setminus \bar \Omega$. Hence $W \big|_{\Omega} \in H^1_0(\Omega)$ is a solution of  $\dive (A \nabla W) + k^2 \Sigma W = 0$ in $\Omega$. We derive that $f = - A \nabla W \cdot \nu \big|_{\Omega} $ on $\partial \Omega$.  This contradicts the fact that $f \not \in {\cal H}$. The proof is complete.  \proofend

\providecommand{\bysame}{\leavevmode\hbox to3em{\hrulefill}\thinspace}
\providecommand{\MR}{\relax\ifhmode\unskip\space\fi MR }
\providecommand{\MRhref}[2]{%
  \href{http://www.ams.org/mathscinet-getitem?mr=#1}{#2}
}
\providecommand{\href}[2]{#2}

\end{document}